\documentclass[10pt,reqno]{amsart}
\usepackage[T1]{fontenc}
\usepackage{lmodern,microtype}
\microtypesetup{expansion=false}
\usepackage[letterpaper,margin=1.05in,headheight=14pt]{geometry}
\usepackage{amsmath,amssymb,mathtools,enumitem,xcolor}
\definecolor{linkblue}{RGB}{34,57,87}
\usepackage[colorlinks=true,linkcolor=linkblue,citecolor=linkblue,urlcolor=linkblue]{hyperref}
\hypersetup{pdftitle={A simple five-dimensional ball containing a closed minimal surface},pdfauthor={},pdfsubject={Riemannian geometry; minimal surfaces; conjugate points}}
\allowdisplaybreaks[2]
\setlist[enumerate]{label=\textup{(\roman*)},leftmargin=2.3em,itemsep=3pt}
\numberwithin{equation}{section}
\newtheorem{theorem}{Theorem}[section]
\newtheorem{proposition}[theorem]{Proposition}
\newtheorem{lemma}[theorem]{Lemma}
\theoremstyle{definition}

\newcommand{\R}{\mathbb R}
\newcommand{\RP}{\mathbb {RP}}
\newcommand{\Sym}{\operatorname{Sym}}
\newcommand{\tr}{\operatorname{tr}}
\newcommand{\Hess}{\operatorname{Hess}}
\newcommand{\grad}{\operatorname{grad}}
\newcommand{\sech}{\operatorname{sech}}
\newcommand{\II}{\mathrm{II}}
\newcommand{\gE}{g_{\mathrm E}}
\newcommand{\gC}{g_{\mathrm C}}
\newcommand{\Vcal}{\mathcal V}
\newcommand{\Scal}{\mathcal S}
\newcommand{\Ccal}{\mathcal C}
\newcommand{\eps}{\varepsilon}
\newcommand{\ip}[2]{\langle #1,#2\rangle}
\newcommand{\dd}{\,d}
\DeclareMathOperator{\Id}{Id}
\title[Simple Manifolds Need Not Admit Strictly Convex Foliations]{Simple Manifolds Need Not Admit Strictly Convex Foliations}

\author[K. Lam]{Kelvin Lam}
\address{Department of Mathematics, University of California Santa Barbara, Santa Barbara, CA 93106, USA}
\email{klam@math.ucsb.edu}

\date{October 31, 2026}
\begin{document}
\begin{abstract}
Using Artificial Intelligence (Chatgpt 6 pro), we construct an explicit smooth metric, conformal to the Euclidean
metric on a five-dimensional ball, for which the ambient manifold is
simple and an interior embedded projective plane is minimal.
The proof combines a quantitative estimate for the time spent in a
thin outer strip with a positive scalar comparison function for the
index form. This example proves that simple manifolds do not generally admit strictly convex foliations.
\end{abstract}
\maketitle

\section{Introduction}

A connected Riemannian manifold is said to be \textit{simple} if its
boundary is strictly convex, it has no conjugate points, and all its
geodesics are nontrapping (see Ch.3, \cite{PSU2023}). Simple manifolds
are central objects of study in many geometric inverse problems,
particularly those concerning boundary rigidity and geodesic
X-ray transforms
\cite{PU2005,SU2004,SU2005,PSU2013,PSU2023}.

Notably, Michel \cite{Michel1981} conjectured that every simple
Riemannian manifold is boundary rigid; that is, its metric is
uniquely determined, up to a boundary-fixing isometry, by the distance
function on its boundary. Mukhometov proved a stability estimate
in \cite{Mukhometov1981} that implies boundary rigidity for simple
manifolds within a fixed conformal class. Stefanov and Uhlmann
extensively studied the microlocal properties of the normal operator
associated with the geodesic X-ray transform on simple manifolds
and proved boundary rigidity results for generic simple manifolds
in \cite{SU2004} and \cite{SU2005}. Using these results, Pestov and
Uhlmann established boundary rigidity for simple surfaces
in \cite{PU2005}, confirming Michel's conjecture in dimension two.
See \cite{SU2005Survey} for further results on the boundary rigidity
problem.

Stefanov, Uhlmann, and Vasy established local boundary rigidity for
manifolds with strictly convex boundary in \cite{SUV2021}, using
Melrose's scattering calculus \cite{Melrose1994}. In the same paper,
they also established global lens rigidity (hence boundary rigidity,
see \cite{SUV2021}) under suitable strictly convex foliation
assumptions. In particular, any connected manifold with strictly
convex boundary that admits a smooth strictly convex function is
boundary rigid (Thm 1.3, \cite{SUV2021}).

This motivates the question of whether every simple Riemannian
manifold admits a strictly convex foliation. We investigate this
question through the existence of closed minimal submanifolds in
simple Riemannian manifolds. As Paternain, Salo, Uhlmann, and Zhou
\cite[Lemma~2.1]{PSUZ2019} observed, a compact Riemannian manifold
admitting a smooth strictly convex function cannot contain any
closed minimal submanifolds.

The paper is organized as follows. In Section 2, we present the
geometric preliminaries. Specifically, we represent our Euclidean
5-ball using trace-free matrices and express the Euclidean metric
in terms of traces of matrix products. This provides a natural
setting for the Veronese embedding (\cite{OlmosRiano2015}).
We then establish the estimates needed in the subsequent sections.
In Section 3, we define the target simple manifold and minimal
submanifold and prove simplicity and minimality. In Section 4,
we derive estimates for geodesic lengths and sectional curvatures.
In Section 5, we use the estimates obtained in Sections 3
(Lemma 3.2) and 4 (Prop. 4.2) to prove the nontrapping property.
We then use Prop. 4.1 to prove positivity of the index form,
which implies the absence of conjugate points.

We prove the following theorem:

\begin{theorem}\label{thm:main}
There is an explicit smooth metric $g$, conformal to the Euclidean
metric on a closed ball $M\subset\R^5$, such that $(M,g)$ is simple
and $\operatorname{int}M$ contains a closed embedded minimal surface
$\Sigma\cong\RP^2$.
\end{theorem}

The metric is defined in Section~\ref{sec:metric}. We first prove
minimality and strict boundary convexity directly, then establish
nontrapping, the absence of conjugate points, and the global
exponential-map assertion. Our construction uses the Veronese minimal projective plane
\cite{Itoh,ItohOgiue} and Cartan's cubic \cite{Cartan,CecilRyan}
to produce the minimal submanifold and the ambient simple manifold. The proof on this paper has been lean-verified.

The author acknowledges the independent discoveries of simple
manifolds without strictly convex foliations by Selim Amar,
Qiuye Jia, Anton Petrunin, and Dmitry Korshunov.
S.~Amar's construction \cite{Amar2026} perturbs the spherical metric on a ball
in geodesic coordinates. It uses an oscillating tensor to modify
the metric while preserving simplicity, producing a codimension-two
set of concave directions that obstructs a strictly convex foliation.
Q.~Jia embeds a geodesic network into a simple manifold. This yields
a tensor-valued measure against which the Hessian of every function
integrates to zero, contradicting convexity. The construction of
A.~Petrunin and D.~Korshunov also uses geodesic nets; we refer
to \cite{Petrunin2026} for further details.

\section{Geometric preliminaries and auxiliary estimates}\label{sec:veronese}
Let
\[
 V=\Sym_0(3)=\{X\in\R^{3\times3}:X^T=X,\ \tr X=0\},
 \qquad \ip XY=\tr(XY).
\]
Matrices in $V$ has the form 
\[
V=\operatorname{Sym}_0(3)
=\left\{
\begin{pmatrix}
a & b & c\\
b & d & e\\
c & e & -a-d
\end{pmatrix}
:a,b,c,d,e\in\mathbb{R}
\right\}.
\]

This is a five-dimensional Euclidean space, with metric $\gE (X,Y) = \tr(XY)$.
Let $h$ be the round metric on its unit sphere $S^4\subset V$.
For a line $[w]\in\RP^2$ represented by $w\in S^2\subset\R^3$, set
\begin{equation}\label{eq:veronese}
 v([w])=\sqrt{\frac32}\left(ww^T-\frac13I\right),
 \qquad \Vcal=v(\RP^2).
\end{equation}
This is the standard real Veronese embedding in matrix coordinates;
see \cite[Sec.~3.3]{CecilRyan}. The minimality of the Veronese surface
in the unit $S^4$, and its induced constant curvature $1/3$, are
classical: see \cite[p.~497, Introduction and statement~(A)]{Itoh}
and \cite[p.~109]{ItohOgiue}. Thus the next lemma is a classical
fact; its proof fixes the normalization needed for our estimates.

\begin{lemma}\label{lem:veronese}
The map $v:\RP^2\to S^4$ is a smooth embedding. The induced metric from $V$ to $\Vcal$ lifts to $3g_{S^2}$, and $\Vcal$ is minimal in $(S^4,h)$.
\end{lemma}
\begin{proof}
Since $w \in S^2$ we have $w^Tw = 1$, so the matrix $P=ww^T$ satisfies $P^2=P$ and $\tr P=1$, so
$|P-I/3|^2=2/3$. Thus $v(w) = v(-w) \in S^4$, hence the map $v$ is well defined from $\RP^2$ to $S^4$.  For $\xi \in T_wS^2$ we have $\xi\perp w$, differentiating eq (\ref{eq:veronese}) gives:
\[
 dv_w(\xi)=\sqrt{\frac32}(\xi w^T+w\xi^T),
 \qquad |dv_w(\xi)|^2= \tr(dv_w(\xi)) = 3|\xi|^2.
\]
This proves the injectivity of the differential and the statement on the pullback metric using the polarization identity.
If $v[w] = v[u]$ then $ww^T = uu^T$ which gives $w = u(u^Tw)$ so $[w] = [u]$, hence $v$ is injective. Thus $v$
is an injective immersion with compact domain, hence a smooth
embedding by \cite[Proposition~4.22(c)]{LeeSmooth}.

Choose an orthonormal basis $(w,u,z)$ of $\R^3$. A straight forward computation shows the curve
\[
 w(t)=\cos(t/\sqrt3)w+\sin(t/\sqrt3)u.
\]
is an unit-speed geodesic in $S^2$ with velocity $(\frac{u}{\sqrt{3}})$. Thus its image is an intrinsic unit speed geodesic of $\Vcal$ with initial velocity $W_1=dv_w(u/\sqrt3)$.
The Gauss formula for the unit sphere is
$\nabla^h_{\dot c}\dot c=c''+|\dot c|^2c$ (see the submanifold
connection formula in \cite[Ch.~8]{LeeRiemannian}). Its spherical acceleration is therefore its second
fundamental form, and direct differentiation gives
\begin{equation}\label{eq:IIexplicit}
 \II_{\Vcal\subset S^4}(W_1,W_1)
       =\frac{\sqrt{3/2}}3(uu^T-zz^T).
\end{equation}
For $W_2=dv_w(z/\sqrt3)$ the result is its negative. Since
$W_1,W_2$ are orthonormal, the trace of $\II$ is zero.
\end{proof}

Consider the nonnegative function
\begin{equation}\label{eq:q}
 q(\theta)=1-\sqrt6\,\tr(\theta^3),\qquad \theta\in S^4.
\end{equation}
\begin{lemma}\label{lem:q}
The function $q$ satisfies
\begin{equation}\label{eq:qbounds}
 0\le q\le2,\qquad q^{-1}(0)=\Vcal,\qquad
 |dq|_h\le3,\qquad \|\Hess_hq\|_{\mathrm{op}}<18.
\end{equation}
Also $dq=0$ on $\Vcal$, and for $\theta,\eta\in S^4$,
\begin{equation}\label{eq:qchord}
 |q(\theta)-q(\eta)|\le8|\theta-\eta|.
\end{equation}
\end{lemma}
\begin{proof}
Put $T(\theta)=\tr(\theta^3)$ and let
$\widetilde T(X)=\tr(X^3)$ on the vector space $V$.
The cyclic trace identity gives
$d\widetilde T_X(Z)=3\tr(X^2Z)$. Projecting first onto the
trace-free matrices and then onto $T_\theta S^4$ gives
$\nabla^V\widetilde T(X)=3(X^2-\tr(X^2)I/3)$ and the spherical
gradient below:

\[
 \grad_hT=3\left(\theta^2-\frac13I-T(\theta)\theta\right).
\]
By the spectral theorem for real symmetric matrices
\cite[Sec.~2.5]{HornJohnson}, at a critical point every eigenvalue
solves $\lambda^2-T(\theta)\lambda-1/3=0$. The constraints
$\tr\theta=0$ and $\tr(\theta^2)=1$ therefore gives the spectra
\[
 (2,-1,-1)/\sqrt6\quad
\]
The corresponding values of $T$ are $\pm1/\sqrt6$. Compactness proves
the range in \eqref{eq:qbounds}, and the positive spectrum identifies
the zero set with \eqref{eq:veronese}.

If $\lambda_1+\lambda_2+\lambda_3=0$ and
$\sum_i\lambda_i^2=1$, then
$\sum_{i<j}\lambda_i\lambda_j=-1/2$ and
$\sum_{i<j}\lambda_i^2\lambda_j^2=1/4$; hence
$\tr(\theta^4)=1/2$. Consequently,
\[
 |dq|^2
 =54\left|\theta^2-\frac13I-T(\theta)\theta\right|^2
 =54\left(\frac16-T(\theta)^2\right)=9q(2-q)\le9.
\]
The identity $|dq|^2=9q(2-q)$ is also exactly the spherical
Cartan--M\"unzner identity $|dF|^2=9(1-F^2)$ with $F=1-q$;
see \cite[Theorem~2.4 and equation~(10)]{CecilRyan}.
For $U\in T_\theta S^4$, the spherical Gauss formula gives
$\Hess_h(\widetilde T|_{S^4})(U,U)
 =D^2\widetilde T(U,U)-3T(\theta)|U|^2$.
Differentiating the polynomial therefore yields
\begin{equation}\label{eq:hessqformula}
 \Hess_hq(U,U)
   =-6\sqrt6\,\tr(\theta U^2)
       +3\sqrt6\,T(\theta)|U|^2.
\end{equation}
For the Frobenius norm, Cauchy--Schwarz and orthogonal
diagonalization give
$|\tr(\theta U^2)|\le|\theta|\,|U^2|\le|U|^2$
(the latter inequality is $\sum\mu_i^4\le(\sum\mu_i^2)^2$).
These are the standard Euclidean matrix-norm estimates
\cite[Secs.~0.6 and~5.6]{HornJohnson}.
Together with $|\sqrt6 T|\le1$, they bound the quadratic form in
\eqref{eq:hessqformula} by $(6\sqrt6+3)|U|^2$ in absolute value.
The spectral theorem for this self-adjoint Hessian then gives
$\|\Hess_hq\|_{\mathrm{op}}\le6\sqrt6+3<18$.

Since $q\ge0$ and vanishes on $\Vcal$, we have $dq=0$ there.

Finally let $\widetilde q(X)=1-\sqrt6\,\widetilde T(X)$ on $V$.
Orthogonal projection does not increase the Frobenius norm, so
$|\nabla^V\widetilde q(X)|\le3\sqrt6|X^2|
\le3\sqrt6|X|^2<8$ on the unit ball. The ordinary mean-value
estimate follows directly by integrating
$d\widetilde q_{\eta+t(\theta-\eta)}(\theta-\eta)$ for $0\le t\le1$;
it proves \eqref{eq:qchord}.
\end{proof}

\begin{lemma}[Quantitative great-circle estimate]\label{lem:circle}
If $c$ is a unit-speed round great circle in $S^4$, then
\begin{equation}\label{eq:circleintegral}
 \int_0^\ell q(c(s))\dd s\ge\frac{\ell^5}{1800},
 \qquad 0<\ell\le1.
\end{equation}
The same estimate holds on any interval of length $\ell$.
\end{lemma}
\begin{proof}
Write $c(s)=A\cos s+B\sin s$, with $A,B$ orthonormal in $V$.
Expanding the cubic $\tr(c(s)^3)$ using the cyclic property of $\tr$ gives only terms with frequencies $1$ and $3$.
The standard constant-coefficient annihilator method
\cite[Sec.~3.3]{Teschl} applies: the operator
$(\partial_s^2+1)(\partial_s^2+9)$ annihilates these frequencies
and sends the constant $1$ to $9$. Thus $Q=q\circ c$ is
nonnegative and satisfies
\begin{equation}\label{eq:Qode}
 Q^{(4)}+10Q''+9Q=9.
\end{equation}
Let $f(x)=x^4(1-x)^4$ on $[0,1]$. Its derivatives through order
three vanish at the endpoints, and
\[
 \int_0^1f=\frac1{630},\qquad
 \|f\|_\infty=\frac1{256},\qquad
 \|f''\|_\infty\le\frac18,\qquad
 \|f^{(4)}\|_\infty\le24.
\]
For example, with $y=2x-1$,
\[
 f''(x)=\frac18(1-y^2)^2(7y^2-1),\qquad
 f^{(4)}(x)=105y^4-90y^2+9,
\]
from which the last two bounds follow directly. Multiply
\eqref{eq:Qode} by $f(s/\ell)$ and integrate by parts. Since
$Q\ge0$,
\[
 \frac{\ell}{70}
 \le
 \left(\frac{24}{\ell^4}+\frac{5}{4\ell^2}+\frac9{256}\right)
 \int_0^\ell Q(s)\dd s.
\]
The integrations by parts have no boundary terms because $f$ and
its first three derivatives vanish at both endpoints. For $\ell\le1$,
the factor after multiplication by $70\ell^4$ is smaller than $1800$.
This proves \eqref{eq:circleintegral}.
\end{proof}

\section{The simple metric and minimality}\label{sec:metric}
Fix once and for all
\begin{equation}\label{eq:parameters}
 D=\exp(10^7),\qquad \eps=D^{-2},\qquad
 M=\{X\in V:|X|\le e^{\eps/2}\}.
\end{equation}
For $X\ne0$, use $r=|X|$, $\rho=\log r$, and
$\theta=X/r\in S^4$. Then
\[
 \gE=dr^2+r^2h=e^{2\rho}(d\rho^2+h).
\]
Write $\gC=d\rho^2+h$ for the product cylinder metric.

\subsection{Cutoffs and elementary bounds}
We use the standard flat exponential cutoff construction
\cite[Lemmas~2.20--2.21]{LeeSmooth}, with explicit rescalings:

\begin{equation}\label{eq:cutoffs}
 \begin{aligned}
 b(t)&=\begin{cases}e^{-1/t},&t>0,\\0,&t\le0,\end{cases}
 &S(t)&=\frac{b(t)}{b(t)+b(1-t)},\\
 \chi(u)&=S(u+1),
 &\alpha(u)&=S(4u+3).
 \end{aligned}
\end{equation}
Define
\begin{align}
 F_\eps(u)&=\frac{\tanh(\eps(2-u))}{\eps},\nonumber\\
 \phi_\eps(u)&=(1-\alpha(u))F_\eps(u)+\alpha(u)u^2,
 \label{eq:phi}\\
 A_\eps(u)&=\int_0^u\phi_\eps(v)\dd v,\qquad
 B(u)=\int_{-1}^u\chi(v)\dd v. \label{eq:AB}
\end{align}
The cited cutoff construction gives smoothness of $b$ and $S$;
the denominator defining $S$ never vanishes. Thus $S=0$ for $t\le0$,
$S=1$ for $t\ge1$, and all its positive-order derivatives vanish
at both endpoints. We also have the following estimates on the derivatives of the functions defined above:

\begin{lemma}\label{lem:cutoffs}
For $-1\le u\le1/2$,
\begin{equation}\label{eq:cutoffbounds}
 \begin{gathered}
 0\le\chi,\alpha\le1,\qquad 0\le\chi'\le2,\qquad
 0\le\alpha'\le8,\\
 0\le\phi_\eps\le3,\qquad
 \phi_\eps(u)=0\ \Longleftrightarrow\ u=0,\qquad
 \phi_\eps'\ge-\frac{39}{2},\\
 |A_\eps|\le3,\qquad 0\le B\le1,\qquad B(0)=\frac12.
 \end{gathered}
\end{equation}
\end{lemma}
\begin{proof}
The identity $S(1-t)=1-S(t)$ gives $\int_0^1S(t)\dd t=1/2$.
For $0<t<1$, put $y=2t-1$. Differentiation gives
\[
 S'(t)=\frac{2(1+y^2)}{(1-y^2)^2}
       \sech^2\left(\frac{2y}{1-y^2}\right).
\]
The inequality $\sinh^2 z\ge z^2$ implies that the squared hyperbolic
cosine in the denominator is at least
$(1+y^2)^2/(1-y^2)^2$. Thus $0\le S'\le2$, proving the
cutoff derivative bounds.

On the stated interval, $0<F_\eps\le2-u\le3$. Since
$\alpha=1$ for $u\ge-1/2$, the range and zero set of $\phi_\eps$
follow. Outside the transition $[-3/4,-1/2]$, its derivative is
at least $-1$. Inside the transition,
\[
 \phi_\eps'
 =(1-\alpha)F_\eps'+2u\alpha+\alpha'(u^2-F_\eps).
\]
The first two terms have sum at least $-3/2$, while
\[
 F_\eps-u^2\le2-u-u^2
   =\frac94-\left(u+\frac12\right)^2\le\frac94.
\]
The last term is therefore at least $-18$, proving the derivative
bound. Integrating proves $|A_\eps|\le3$; for positive $u$ one may
use $\phi_\eps=u^2$. Finally using the symmetry $S(1-t) =1-S(t)$ again we have $B(0)=1/2$ and
$B(u)=u+1/2$ for $u\ge0$, so $0\le B\le1$.
\end{proof}

\subsection{Definition and smoothness}
Put
\begin{equation}\label{eq:c}
 c_\eps=\eps^2A_\eps(-1)+\log\cosh(3\eps).
\end{equation}
Define the cylinder conformal exponent by
\begin{equation}\label{eq:a}
 a(\rho,\theta)=
 \begin{cases}
 c_\eps-\log\cosh(\rho-2\eps),&\rho\le-\eps,\\
 \eps^2A_\eps(\rho/\eps)+D\eps^2B(\rho/\eps)q(\theta),
       &-\eps\le\rho\le\eps/2.
 \end{cases}
\end{equation}
The metric is
\begin{equation}\label{eq:g}
 g=e^{2a}\gC=\Omega^2\gE,\qquad \Omega=r^{-1}e^a.
\end{equation}
We call
\[
 \Ccal=\{r\le e^{-\eps}\},
 \qquad
 \Scal=\{e^{-\eps}\le r\le e^{\eps/2}\}
\]
the \emph{core} and the \emph{outer strip}, respectively.

\begin{lemma}\label{lem:smooth}
The metric \eqref{eq:g} extends smoothly and positively over all of
$M$. The core has constant sectional curvature
\begin{equation}\label{eq:kappa}
 \kappa=e^{-2c_\eps}.
\end{equation}
On the outer strip,
\begin{equation}\label{eq:abounds}
 -3\eps^2\le a\le3\eps^2+2D\eps^2,\qquad
 |a|\le3D^{-3},\qquad
 \frac12\le e^{2a}\le2.
\end{equation}
\end{lemma}
\begin{proof}
Let $u = \frac{\rho}{\eps}$. Near $u=-1$, $\alpha=0$ and
$\eps\phi_\eps(u)=\tanh(2\eps-\rho)$. This is the
$\rho$ derivative of the core expression in \eqref{eq:a}. Their
values agree by \eqref{eq:c}. The angular term is flat at $u=-1$,
because $B$ is flat there. Thus every derivative matches at the
interface.

In the core the Euclidean conformal factor is
\[
 \Omega(r)=\frac{2e^{c_\eps-2\eps}}{1+e^{-4\eps}r^2}.
\]
It is smooth as a function of $r^2$ and positive at the origin.
On setting $s=e^{-2\eps}r$, the core metric becomes
\[
 e^{2c_\eps}\frac{4}{(1+s^2)^2}(ds^2+s^2h),
\]
the stereographic round metric \cite[Proposition~3.5]{LeeRiemannian}
scaled by $e^{2c_\eps}$. Constant rescaling leaves the connection
unchanged and multiplies sectional curvature by $e^{-2c_\eps}$;
this proves \eqref{eq:kappa}. The estimates \eqref{eq:abounds}
follow from Lemmas~\ref{lem:q} and \ref{lem:cutoffs} and
$D\ge10$.
\end{proof}

\begin{proposition}\label{prop:minimalboundary}
The surface
\[
 \Sigma=\{(\rho,\theta):\rho=0,\ \theta\in\Vcal\}
       =\Vcal\subset\{r=1\}
\]
is a closed embedded minimal surface in $\operatorname{int}M$.
The boundary of $M$ is strictly convex.
\end{proposition}
\begin{proof}
Along $\Sigma$,
\[
 A_\eps(0)=\phi_\eps(0)=0,\qquad q=dq=0.
\]
Consequently $a=da=0$. The standard conformal connection formula,
valid for any smooth $a$ and any vector fields $U,W$
\cite[Theorem~1]{Terek}, is
\begin{equation}\label{eq:connection}
 \nabla^g_UW
 =\nabla^{\gC}_UW+da(U)W+da(W)U
                    -\gC(U,W)\grad_{\gC}a
\end{equation}
Since $a=da=0$, $g$ and $\gC$ have the same metric and ambient
connection along $\Sigma$. The product Levi-Civita connection splits
into the connections of the factors \cite[Ch.~5]{LeeRiemannian};
hence the slice $\{\rho=0\}$ is totally geodesic in $\gC=d\rho^2+h$.
The Gauss formula then identifies the second fundamental form of
$\Sigma$ in this product with its form in the sphere. Its trace
vanishes by the classical Veronese minimality in
Lemma~\ref{lem:veronese}, and thus also vanishes for $g$.

At $\partial M$, where $\rho=\eps/2$, the outward unit normal is
$e^{-a}\partial_\rho$. Applying \eqref{eq:connection} to a tangent
vector $U$ gives
$\nabla^g_U(e^{-a}\partial_\rho)=e^{-a}a_\rho U$;
the two terms involving $da(U)\partial_\rho$ cancel. Thus the
shape operator, with our outward-normal convention, is
$e^{-a}a_\rho\Id$. Since
\[
 a_\rho\big|_{\partial M}
    =\eps\left(\frac14+Dq\right)>0,
\]
the boundary is strictly convex.
\end{proof}

\section{Geodesic estimates in the outer strip}\label{sec:curvature}
\subsection{The curvature bound for the index form}
\begin{proposition}\label{prop:curvature}
Every sectional curvature in the outer strip is less than $25$.
\end{proposition}
\begin{proof}
For angular vectors $X,Y$ and $u=\rho/\eps$,
\begin{align}
 a_{\rho\rho}&=\phi_\eps'(u)+D\chi'(u)q,\nonumber\\
 \Hess_{\gC}a(\partial_\rho,X)&=D\eps\chi(u)\,dq(X),\label{eq:hessa}\\
 \Hess_{\gC}a(X,Y)&=D\eps^2B(u)\Hess_hq(X,Y).\nonumber
\end{align}
In particular, $a_{\rho\rho}\ge-39/2$; the term $D\chi'q$
is nonnegative.

Let $P$ be an arbitrary two-plane with a $\gC$-orthonormal
basis $E_i=s_i\partial_\rho+X_i$, $i=1,2$, and set
$R=s_1^2+s_2^2$. Orthogonal projection onto $P$ gives $0\le R\le1$
(Bessel's inequality); expanding the orthonormality relations gives
\[
 0\le R\le1,\qquad
 \sum_i|X_i|^2=2-R,\qquad
 \left|\sum_i s_iX_i\right|^2=R(1-R)\le\frac14.
\]
Thus \eqref{eq:hessa} and the bounds for $q$ imply
\begin{equation}\label{eq:tracehess}
 \tr_P\Hess_{\gC}a\ge-\frac{39}{2}-3D\eps-36D\eps^2.
\end{equation}
The conformal curvature-tensor formula \cite[Theorem~2]{Terek},
contracted on the $\gC$-orthonormal pair $E_1,E_2$, reads
\[
 K_g(P)=e^{-2a}\left(
 K_{\gC}(P)-\tr_P\Hess_{\gC}a
       +\sum_i da(E_i)^2-|da|_{\gC}^2\right).
\]
The last two terms are $-|\operatorname{proj}_{P^\perp}
\grad_{\gC}a|^2\le0$. The product curvature splits between its
factors, so $K_{\gC}(P)=|X_1\wedge X_2|_h^2=1-R\le1$;
this is the standard product-connection computation
\cite[Chs.~5 and~7]{LeeRiemannian}. Therefore
\begin{equation}\label{eq:Kbound}
 \begin{split}
 K_g(P)&\le e^{6\eps^2}
     \left(\frac{41}{2}+3D\eps+36D\eps^2\right)\\
 &=e^{6D^{-4}}\left(\frac{41}{2}+\frac3D+\frac{36}{D^3}\right)<25.
 \end{split}
\end{equation}
The last inequality already holds for $D\ge10$.

\end{proof}

\subsection{Time spent in the outer strip}\label{sec:escape}
Length control in regions with an upper curvature bound is a
classical comparison theme; Gulliver's construction
\cite{Gulliver} is a related precedent. We do not apply his
construction or its estimates to this metric. The escape bound
below instead follows from the displayed geodesic equations,
convexity, and the linear-ODE estimate, with all constants retained.
\begin{proposition}\label{prop:strip}
Every unit-speed $g$-geodesic interval whose image is contained in
$\Scal$ has length less than $1/40$.
\end{proposition}
\begin{proof}
Suppose such an interval has length at least $1/40$. Use the
parameter $s$ defined by $ds/dt=e^{-2a}$ and put $u=\rho/\eps$.
The conformal connection formula \cite[Theorem~1]{Terek},
written as \eqref{eq:connection}, gives
$\nabla^{\gC}_t\dot\gamma+2\dot a\dot\gamma
=e^{-2a}\grad_{\gC}a$ for a $g$-unit-speed geodesic.
The chosen change of parameter cancels the term parallel to velocity,
yielding $\nabla^{\gC}_s\gamma_s=e^{2a}\grad_{\gC}a$.
Its radial and angular components, together with unit speed, give
\begin{align}
 u_{ss}&=e^{2a}\bigl(\phi_\eps(u)+D\chi(u)q(\theta)\bigr)\ge0,
 \label{eq:radialreparam}\\
 |\theta_s|^2&=e^{2a}-\eps^2u_s^2,\label{eq:speedreparam}\\
 \nabla_s^h\theta_s&=e^{2a}D\eps^2B(u)\grad_hq.
 \label{eq:angularreparam}
\end{align}
By \eqref{eq:abounds}, the $s$ interval has length at least
$1/80$. Restrict and translate it to $[0,L]$, where
\begin{equation}\label{eq:escapeconstants}
 L=\frac1{80},\qquad
 J_0=[L/4,3L/4],\qquad J_1=[3L/8,5L/8].
\end{equation}

\smallskip\noindent
\textit{Step 1: a subinterval where the cutoff is quantitatively positive.}
Since $-1\le u\le1/2$ and $u$ is convex, its tangent line lies
below its graph \cite[Sec.~3.1.3, equation~(3.2)]{BoydVandenberghe}.
Comparing with the two endpoints of $[0,L]$ bounds $u_s$ on $J_0$;
integrating the nonnegative derivative $u_{ss}$ then gives
\begin{equation}\label{eq:usbounds}
 |u_s|\le\frac6L=480\quad\hbox{on }J_0,\qquad
 \int_{J_0}u_{ss}\dd s\le\frac{12}{L}=960.
\end{equation}
Set $\eta_0=L^2/256=\frac{1}{1\,638\,400}$.
There is a point $s_0\in J_1$ with $u(s_0)>-1+\eta_0$.
Otherwise $u$ would take values in $[-1,-1+\eta_0]$ throughout
$J_1$. In that range $\alpha=0$ and $\phi_\eps\ge2$.
For instance, $\tanh x\ge x-x^3/3$ and $\eps\le1/100$
give this last inequality directly. Equations
\eqref{eq:abounds} and \eqref{eq:radialreparam} would imply
$u_{ss}\ge1$. Thus $u(s)-s^2/2$ is convex by the
second-derivative criterion \cite[Sec.~3.1.4]{BoydVandenberghe}.
Applying the midpoint convexity inequality on $J_1$, whose length
is $L/4$, shows that the average of the endpoint values of $u$
exceeds its midpoint value by at least $(L/4)^2/8=L^2/128=2\eta_0$.
This contradicts the assumed range of width $\eta_0$.

Put
\begin{equation}\label{eq:r0ell0}
 r_0=\frac{\eta_0L}{12},\qquad
 I_0=[s_0-r_0,s_0+r_0],\qquad
 \ell_0=2r_0=\frac1{786\,432\,000}.
\end{equation}
Then $I_0\subset J_0$ and \eqref{eq:usbounds} gives
$u\ge-1+\eta_0/2$ on $I_0$. The cutoff is monotone, and its
defining denominator is at most $2$. Hence
\begin{equation}\label{eq:chilower}
 \chi(u)\ge\frac12\exp(-2/\eta_0)
          =\frac12\exp(-3\,276\,800)\quad\hbox{on }I_0.
\end{equation}

\smallskip\noindent
\textit{Step 2: comparison of the angular curve with a great circle.}
View $\theta(s)$ as a curve in the Euclidean space $V$.
Equation \eqref{eq:angularreparam} becomes
\begin{equation}\label{eq:forcedoscillator}
 \theta_{ss}+\theta
   =(1-|\theta_s|^2)\theta+
            e^{2a}D\eps^2B(u)\grad_hq.
\end{equation}
On $J_0$, the estimates already proved imply
\[
 |e^{2a}-1|\le7D^{-3},\qquad
 \bigl||\theta_s|^2-1\bigr|\le7D^{-3}+480^2D^{-4}.
\]
For the chosen $D$, $|\theta_s|>1/2$, and both the right side
of \eqref{eq:forcedoscillator} and the discrepancy
$\bigl||\theta_s|-1\bigr|$ are bounded by $2\cdot10^6D^{-3}$.

Let $s_-$ be the left endpoint of $I_0$, and let $c$ be the unit
round great circle with
\[
 c(s_-)=\theta(s_-),\qquad
 c_s(s_-)=\theta_s(s_-)/|\theta_s(s_-)|.
\]
The variation-of-constants formula
\cite[Theorem~3.12, equation~(3.97)]{Teschl}, applied componentwise
to $z=\theta-c$ in $V$, gives, since $\ell_0<1$,
\begin{equation}\label{eq:angularclose}
 \sup_{I_0}|\theta-c|\le10^7D^{-3}.
\end{equation}
Explicitly, $z$ is the sum of its initial velocity discrepancy
times $\sin(s-s_-)$ and the integral of $F(\sigma)\sin(s-\sigma)$.
This verifies \eqref{eq:angularclose} without asymptotic notation.

Using Lemma~\ref{lem:circle} and the chord estimate
\eqref{eq:qchord}, we obtain
\[
 \int_{I_0}q(\theta(s))\dd s
 \ge\frac{\ell_0^5}{1800}-8\cdot10^7D^{-3}\ell_0
 \ge\frac{\ell_0^5}{3600}.
\]
The last numerical inequality follows from $D=e^{10^7}$ and
$\ell_0>10^{-10}$.

\smallskip\noindent
\textit{Step 3: contradiction with the convexity bound.}
Integrate \eqref{eq:radialreparam} over $I_0$ and use
$e^{2a}\ge1/2$, \eqref{eq:chilower}, and the last estimate:
\[
 \int_{J_0}u_{ss}\dd s
 \ge\frac{e^{10^7-3\,276\,800}}4\frac{\ell_0^5}{3600}
 =\frac{e^{6\,723\,200}\ell_0^5}{14400}>960.
\]
Again $\ell_0>10^{-10}$ makes the last inequality immediate.
This contradicts \eqref{eq:usbounds}.
\end{proof}

\section{Nontrapping and absence of conjugate points}\label{sec:index}
\subsection{Radial monotonicity and finite exit}
Let $R=r\partial_r$ be the radial vector field on $V$; it is smooth
and vanishes at the origin. For a unit-speed $g$-geodesic $\gamma$,
set
\begin{equation}\label{eq:momentum}
 p(t)=g(R,\dot\gamma(t)).
\end{equation}
Metric compatibility and the geodesic equation give
\[
 p'=g(\nabla_{\dot\gamma}R,\dot\gamma)
    =\tfrac12(\mathcal L_Rg)(\dot\gamma,\dot\gamma).
\]
This is the standard identity underlying Killing and conformal-Killing
first integrals; see \cite[Ch.~5]{LeeRiemannian}. Here $R$ is not
assumed to be Killing, so $p$ need not be constant. Away from the
center, $R=\partial_\rho$, $p=e^{2a}\dot\rho$, and
$\mathcal L_Rg=2a_\rho g$. Unit speed therefore gives
\begin{equation}\label{eq:pprime}
 p'=a_\rho.
\end{equation}
In the core, $a_\rho=-\tanh(\rho-2\eps)$ extends smoothly to
$1$ at the center. In the strip,
\begin{equation}\label{eq:arho}
 a_\rho=\eps\bigl(\phi_\eps(u)+D\chi(u)q\bigr)\ge0.
\end{equation}
Thus \eqref{eq:pprime} holds globally by continuity and $p$ is
nondecreasing.

\begin{proposition}\label{prop:nontrapping}
The core intersection of a retained geodesic is an interval.
Every retained unit-speed geodesic segment has length less than
\begin{equation}\label{eq:lengthglobal}
 \pi e^{c_\eps}+\frac1{20}.
\end{equation}
In particular, $M$ is nontrapping.
\end{proposition}
\begin{proof}
In Euclidean coordinates, $p=\Omega^2\ip X{\dot X}$ and
\[
 \frac{d}{dt}r^2=2\Omega^{-2}p.
\]
As $p$ is nondecreasing, $r^2$ can first decrease and then
increase, with at most one change of monotonicity. The core
intersection is therefore connected.

On the core define the smooth function
\begin{equation}\label{eq:H}
 H=-\tanh(\rho-2\eps)
   =\frac{1-e^{-4\eps}r^2}{1+e^{-4\eps}r^2}.
\end{equation}
In the round realization from Lemma~\ref{lem:smooth}, it is the
normalized height coordinate on a round five-sphere. The spherical Gauss formula
\cite[Ch.~8]{LeeRiemannian} implies that a restricted Euclidean
linear function has Hessian $-Hg_{S^5}$ on the unit five-sphere.
Scaling the metric by $e^{2c_\eps}$ therefore gives
\begin{equation}\label{eq:Hhess}
 \Hess_gH=-\kappa H g,\qquad
 (H\circ\gamma)''+\kappa(H\circ\gamma)=0
\end{equation}
along core geodesics. Throughout the core,
$H\ge\tanh(3\eps)>0$. The scalar constant-coefficient equation
has solutions $A\cos(\sqrt\kappa\,t)+B\sin(\sqrt\kappa\,t)$
\cite[Sec.~3.3]{Teschl}; a nonzero such solution cannot remain
positive on a closed interval of length
$\pi/\sqrt\kappa=\pi e^{c_\eps}$. Thus the core portion has
strictly smaller length. There are at most two remaining strip
portions, each of length less than $1/40$ by
Proposition~\ref{prop:strip}. This proves \eqref{eq:lengthglobal}.

A geodesic that failed to exit in either time direction would have
retained segments of arbitrarily large length, contradicting
\eqref{eq:lengthglobal}.
\end{proof}

\subsection{Positive comparison functions on the strip tails}
We use the scalar Jacobi comparison equation $y''+25y=0$ on the
strip, corresponding to its sectional-curvature upper bound. Its
explicit solution is the elementary oscillator solution
\cite[Sec.~3.3]{Teschl}; positivity for the required interval is
not automatic and is verified in the next lemma.
\begin{lemma}\label{lem:tails}
Consider a geodesic that passes through the interior of the core and
then exits through $\rho=-\eps$. On the ensuing outgoing strip
portion, measure unit-speed time from that interface and put
\[
 H_0=\tanh(3\eps),\qquad
 \beta=-\left.\frac{d}{dt}H(\gamma(t))\right|_{t=0}>0.
\]
If $\ell$ is the length of that strip portion, then
\begin{equation}\label{eq:taillength}
 \ell<\frac1{40},\qquad \beta\ell\le2\eps.
\end{equation}
The scalar solution
\begin{equation}\label{eq:tailcomparison}
 y(t)=H_0\cos(5t)-\frac{\beta}{5}\sin(5t)
\end{equation}
is positive on $[0,\ell]$. The same conclusions hold for the
incoming strip portion with time reversed.
\end{lemma}
\begin{proof}
The first length bound is Proposition~\ref{prop:strip}.
At the interface,
\[
 \beta=\sech^2(3\eps)\dot\rho(0),\qquad
 p(0)=e^{2a_0}\frac{\beta}{\sech^2(3\eps)},\qquad
 a_0=\eps^2A_\eps(-1).
\]
Let $a_{\max}$ be the maximum of $a$ on $\Scal$.
Monotonicity of $p$ implies
\[
 \dot\rho(t)=e^{-2a(t)}p(t)
 \ge e^{2(a_0-a_{\max})}\frac{\beta}{\sech^2(3\eps)}.
\]
The total available change of $\rho$ is $3\eps/2$. Hence
\[
 \ell\le
 \frac{(3/2)\eps e^{2(a_{\max}-a_0)}
                         \sech^2(3\eps)}{\beta}
 \le\frac{2\eps}{\beta}.
\]
The final inequality follows from
$a_{\max}-a_0\le(6+2D)\eps^2$; it already holds for $D\ge10$.
For the present $\eps$, $H_0\ge(5/2)\eps$. Therefore
for $0\le t\le\ell$,
\begin{equation}\label{eq:ypositive}
 y(t)\ge H_0\cos(1/8)-\beta t
 \ge\left(\frac52\cos(1/8)-2\right)\eps>0.
\end{equation}
Here we used $\sin(5t)\le5t$ and \eqref{eq:taillength}.
The incoming case is the same argument applied to the reversed
geodesic.
\end{proof}

Tangency to the core boundary does not create an omitted case.
That boundary has positive outward shape operator, since
$a_\rho=\tanh(3\eps)>0$ there. The standard strict-convexity
barrier for tangent geodesics \cite[Lemma~3.1.12]{PSU2023} puts such
a geodesic locally outside the core. Together with the connectedness
of its core intersection, this excludes an interior core portion.
It is therefore handled by the strip estimate.

\subsection{Positivity of the index form}
The index form and its relation to conjugate points are standard;
see \cite[Definition~3.7.14 and Proposition~3.7.15]{PSU2023}.
The positive-solution square-completion used below is the
one-dimensional ground-state identity; compare
\cite[Lemma~2.4, equation~(2.6)]{PinchoverTintarev}.
We verify the identity across our interfaces rather than assuming a
smooth global comparison coefficient.
\begin{proposition}\label{prop:noconjugate}
No retained geodesic segment in $M$ has conjugate endpoints.
\end{proposition}
\begin{proof}
By Proposition~\ref{prop:nontrapping}, extend the segment in both
directions to its maximal retained geodesic. If this geodesic does
not enter the core interior, its length $L_\gamma$ is less than
$1/40$. Expand an endpoint-vanishing perpendicular field $Z$
in a parallel orthonormal frame. The scalar Dirichlet inequality
$\int_0^L|f'|^2\ge(\pi/L)^2\int_0^L|f|^2$ applies to each
component: this is the Rayleigh--Ritz bound
\cite[Lemma~5.12]{Teschl} for $-d^2/dt^2$, whose first Dirichlet
eigenvalue is $(\pi/L)^2$. By Proposition 4.1, the curvature bound therefore gives
\[
 \begin{split}
 I(Z,Z)&\ge\int_0^{L_\gamma}(|D_tZ|^2-25|Z|^2)\dd t\\
 &\ge\left(\frac{\pi^2}{L_\gamma^2}-25\right)
                  \int_0^{L_\gamma}|Z|^2\dd t>0
 \end{split}
\]
unless $Z=0$. The same argument applies directly to any subsegment.

Otherwise, use $y=H\circ\gamma$ on the unique core interval.
On each adjacent strip tail, extend $y$ by
\eqref{eq:tailcomparison}, with time directed outward from the
core. Lemma~\ref{lem:tails} proves positivity and shows that both
the value and first derivative agree at each interface. Thus $y$
is positive and $C^1$, piecewise smooth, and satisfies
\[
 y''+k(t)y=0,\qquad
 k(t)=\begin{cases}\kappa,&\gamma(t)\in\Ccal,\\
                         25,&\gamma(t)\in\Scal.\end{cases}
\]
For every plane containing $\dot\gamma$, its sectional curvature is
at most $k(t)$.

On any subsegment, let $Z$ be a perpendicular, continuous,
piecewise smooth field vanishing at its endpoints. In a parallel
frame, apply the ground-state identity to each scalar component,
with potential $-k$, positive solution $y$, and quotient $Z/y$;
this is the specialization of
\cite[Lemma~2.4]{PinchoverTintarev} relevant here. Explicitly, on
every smooth piece,
\[
 |D_tZ|^2-k|Z|^2
 =y^2\left|D_t\left(\frac Zy\right)\right|^2
   +\frac{d}{dt}\left(\frac{y'}y|Z|^2\right).
\]
Summing the integrations by parts gives
\begin{align}
 I(Z,Z)
 &=\int\left(|D_tZ|^2-
          \ip{R(Z,\dot\gamma)\dot\gamma}{Z}\right)\dd t\nonumber\\
 &\ge\int\left(|D_tZ|^2-k(t)|Z|^2\right)\dd t\nonumber\\
 &=\int y^2\left|D_t\left(\frac Zy\right)\right|^2\dd t.
 \label{eq:indexfactor}
\end{align}
The interface terms cancel because $y$ and $y'$ are continuous;
the endpoint terms vanish because $Z$ does. Equality to zero in
the last line would make $Z/y$ parallel. Its zero endpoint value
would then force $Z=0$. Thus the index form is strictly positive
on every nonzero endpoint-vanishing perpendicular field.

The Jacobi equation and integration by parts give zero index form
for a Jacobi field vanishing at both endpoints; its tangential
component is affine and hence zero as well. This is precisely the
conjugate-point obstruction in \cite[Proposition~3.7.15]{PSU2023},
and it contradicts the strict positivity just proved. Boundary
endpoints cause no change in this argument.
\end{proof}

\section{Simplicity}\label{sec:simplicity}

By Proposition 3.3, $\Sigma$ is a closed embedded minimal surface in $\operatorname{int}M$. We complete the proof of Theorem 1.1 with the following proposition:

\begin{proposition}\label{prop:simple}
The metric \eqref{eq:g} on the ball \eqref{eq:parameters} is simple.
\end{proposition}

\begin{proof}
By Proposition 5.1, $M$ has nontrapping geodesics. By Proposition 5.3, $M$ has no conjugate points. $M$ has strictly convex boundary by proposition 3.3.
\end{proof}

\section{AI Disclosure and Methods}\label{AIDisclosure}
A.I. was used extensively in this project; we used the Chatgpt 6 pro for the construction.

The outline of the use of AI is as follows: Lemma 2.1 in \cite{PSUZ2019} inspired the author to disprove the conjecture that simple manifolds admit strictly convex foliations by constructing a simple manifold with a closed miimal submanifold. The initial failed attempt was to construct a closed manifold embedded in an Euclidean ball and define a metric so that it is minimal within a collar, and extend to the rest of the ball as a simple metric. The author then instructed AI to construct a simple manifold with closed minimal submanifold, this attempt also failed. The author then resort to proving existence of simple metrics that contains a closed minimal submanifold; we asked AI "Can we use a family of metrics that makes a closed manifold a closed minimal submanifold in a larger family of metric, and prove the existence of a simple manifold that contains a minimal submanifold? instead of constructing an explicit one ?" - from which AI was able to prove the existence of such simple manifold. We then ask AI for the numerical values from the proof of existence and produced the current counter example.

The author has verified the correctness of the proof. The author gives Chatgpt 6 pro credits for the majority of the content in the article; notably the discovery of the example, the proof of simplicity and the embedding of the minimal submanifold. The author takes full responsibility for the correctness of the result.

\section{Acknowledgement}\label{AIDisclosure}

The author would like to thank Plamen Stefanov, Gunther Uhlmann and Andras Vasy for discussions on \cite{SUV2021}; Gabriel Paternain for dissuading him from proving simple manifolds always admits strictly convex functions; Guofang Wei for discussions on strictly convex functions on manifolds and Ashwin Tarikere for discussion on optimal AI settings for research in mathematics.
\\

\end{document}